\documentclass[12pt]{amsart}

\usepackage{amsmath,amsfonts,amsthm,amsopn,cite,mathrsfs}
\usepackage{epsfig,verbatim}
\usepackage{subfigure}

\usepackage{amssymb}
\usepackage{amsthm}
\newtheorem{thm}{Theorem}[section]

\newtheorem{defi}[thm]{Definition}
\newtheorem{lem}[thm]{Lemma}

\renewcommand{\leq}{\leqslant}

\renewcommand{\geq}{\geqslant}
\newcommand\U{\mathcal{U}}
\newcommand\G{\mathcal{G}}
\newcommand\s{\mathcal{S}}
\newcommand\PP{\mathcal{P}}
\newcommand\N{\mathcal{N}}

\newcommand\R{\mathcal{R}}
\newcommand{\ttt}[1]{\mbox{\textnormal{(\ref{#1})}}}

\begin{document}

\title[Sampling and Reconstruction by Oblique Projections]{Sampling and Reconstruction in Different Subspaces by  Using
  Oblique Projections}

\author{Peter Berger}
\address{Faculty of Mathematics \\
University of Vienna \\
Oskar-Morgenstern-Platz 1 \\
A-1090 Vienna, Austria}
\email{p.berger@univie.ac.at}
\author{Karlheinz Gr\"ochenig}
\address{Faculty of Mathematics \\
University of Vienna \\
Oskar-Morgenstern-Platz 1\\
A-1090 Vienna, Austria}
\email{karlheinz.groechenig@univie.ac.at}
\subjclass[2000]{}
\date{}
\keywords{}
\thanks{The authors were 
  supported by the  National Research Network S106 SISE of the
Austrian Science Fund (FWF)} 
\maketitle

\section{Introduction}
We consider an approach, where sampling and reconstruction are done in different subspaces of a Hilbert space $H$.
Let $U$ be a sampling subspace, and let $G$ be a reconstruction subspace.
Let $\{u_j\}_{j\in \mathbb{N}}$ be a frame for $U$. Given the scalar products 
 $\{\langle f,u_j\rangle\}_{j\in\mathbb{N}}$ of an element $f\in H$, we want to find a stable reconstruction $\tilde{f}$ of $f$,
 that is close to the unknown orthogonal projection of $f$ onto the reconstruction space $G$.
 
The main idea of our reconstruction method is explained in the
following example. Let $H=\mathbb{R}^3$, let $G$ be a one dimensional
subspace of $H$,  
and let $U$ be the linear span of two linearly independent vectors $u_1$ and $u_2$. We intend to reconstruct $f\in H$ from
$\langle f,u_1\rangle$ and $\langle f,u_2\rangle$. From the measurements $\langle f,u_1\rangle$ and $\langle f,u_2\rangle$, 
we can calculate $P_Uf$, the projection of $f$ onto the plane $U$. Conversely, $P_U f$ determines
$\langle f,u_1\rangle$ and $\langle f,u_2\rangle$. 

Thus all the information we have about $f$ is that $f$ lies in the  affine subspace $P_Uf+U^\perp$,
but we do not know the exact location of $f$ in this affine subspace. 
Let $\mathcal{P}_{G P_U(G)^\perp}$ denote the oblique projection with
range $G$ and kernel $P_U(G)^\perp$. 

We assume that $f$, the element to be  reconstructed, is close the reconstruction space $G$.
Naturally, we now want to find $\tilde{f}$, the element of $G$ (the reconstruction space) closest to $P_Uf+U^\perp$.
The two spaces $P_Uf+U^\perp$ and $G$ may, or may not intersect. 
In both cases, the element of $G$ closest to $P_Uf+U^\perp$ is exactly $\mathcal{P}_{G P_U(G)^\perp}f$.
If they intersect, then $\mathcal{P}_{G P_U(G)^\perp}f = (P_Uf+U^\perp)\cap G$, and
$\langle f,u_1\rangle = \langle \tilde{f},u_1\rangle$ and $\langle f,u_2\rangle = \langle \tilde{f},u_2\rangle$.
In this case $\tilde{f}$ is a so called consistent reconstruction of $f$.

One should remember that in this setup only the scalar products of $f$ with the frame sequence 
$\{u_j\}_{j\in\mathbb{N}}$ of $U$ are given. Thus we analyse the operator 
$Q:\{\langle f,u_j\rangle\}_{j\in\mathbb{N}}\mapsto \mathcal{P}_{G P_U(G)^\perp}f$.
 An explicit formula for the mapping $Q$ is given   in Theorem \ref{wichtig}. 
We refer to the mapping $\mathcal{P}_{G P_U(G)^\perp}$ as
\textit{frame independent sampling} to  
indicate that $\mathcal{P}_{G P_U(G)^\perp}$ does not depend on the frame sequences $\{u_j\}_{j\in \mathbb{N}}$ and
$\{g_k\}_{j\in \mathbb{N}}$ themselves, but only on their closed linear spans $U$ and $G$.

This mapping $Q$ is a generalisation of 
consistent sampling, which is treated in \cite{23,24,28}. If the frame sequence $\{u_j\}_{j\in\mathbb{N}}$ of the sampling 
space $U$ is tight, this reconstruction coincides with the generalized sampling introduced in \cite{1,adcock,adhapo12}.
In the following, we study this operator, and compare it with the generalized sampling introduced in 
\cite{1,adcock,adhapo12}.

\section{Stability and quasi-optimality} \label{stabo}
Let $\{u_j\}_{j\in \mathbb{N}}$ be a frame sequence in $H$, i.e., a frame for its closed linear span. Setting  
$U:=\overline{\textnormal{span}}\{u_j\}_{j\in \mathbb{N}}$, this is equivalent to the statement that there exist constants $A,B>0$, such that 
\begin{equation} \label{frameequation}
  A \|f\|^2\leq \sum_{j\in \mathbb{N}}|\langle f,u_j\rangle |^2\leq B\|f\|^2,\quad \text{for every }f\in U.
 \end{equation}
The constant $A$ is called lower frame bound and the constant $B$ is called upper frame bound.   
We call 
$$U=\overline{\textnormal{span}}\{u_j\}_{j\in \mathbb{N}}$$
the \textit{sampling space}. 

Let $\{g_k\}_{k\in \mathbb{N}}$ be a frame sequence in $H$. Setting  
$U:=\overline{\textnormal{span}}\{g_k\}_{k\in \mathbb{N}}$, this is equivalent to the statement that there exist constants $C,D>0$, such that 
\begin{equation}
  C \|{f}\|^2\leq \sum_{k\in \mathbb{N}}|\langle f,g_k\rangle |^2\leq D\|f\|^2,\quad \text{for every }f\in G.
 \end{equation}
We call 
$$G=\overline{\textnormal{span}}\{g_k\}_{j\in \mathbb{N}}$$
the \textit{reconstruction space}.
 The operator 
\begin{equation}
 \mathcal{U}:l^2(\mathbb{N})\rightarrow H,\quad \mathcal{U}\{c_j\}_{j\in \mathbb{N}}=\sum_{j=1}^\infty c_j u_j
\end{equation}
 is called the \textit{synthesis operator} of the frame sequence $\{u_j\}_{j\in \mathbb{N}}$.
The adjoint operator
\begin{equation}
  \mathcal{U}^*:H\rightarrow l^2(\mathbb{N}),\quad \mathcal{U}^*f= \{\langle f,u_j\rangle\}_{j\in \mathbb{N}}
\end{equation}
is called the \textit{analysis operator} of the frame sequence $\{u_j\}_{j\in \mathbb{N}}$.
The composition 
\begin{equation}
\mathcal{S}:H\rightarrow H,\quad \mathcal{S}f=\mathcal{U}\mathcal{U}^*f=\sum_{j=1}^\infty \langle f,u_j\rangle u_j
\end{equation}
is called the frame operator.

From now on we denote by $\U$ the synthesis operator, by $\U^*$ 
the analysis operator, and by $\s$ the frame operator of the frame sequence $\{u_j\}_{j\in \mathbb{N}}$ of the sampling space $U$. 

We denote by $\G$ the synthesis operator and by $\G^*$ 
the analysis operator of the frame sequence $\{g_k\}_{k\in \mathbb{N}}$ of the reconstruction space $G$.

In the following we denote by $\R(A)$ the range of the operator $A$ and by $\N(A)$ the nullspace of the operator $A$.

As mentioned in the introduction, we want  to find a mapping 
\begin{equation}
Q:\R(\U^*)\rightarrow G,
\end{equation}
such that the mapping 
\begin{equation}
F:=Q\U^*
\end{equation}
\begin{equation}
F:H\rightarrow G, \quad Ff=Q\U^*f,
\end{equation}
has the property that $Ff$ is a good approximation to $f$ for every $f\in H$. Following \cite{adhapo12}, we use two quantities to measure the quality of the 
reconstruction $F=Q\U^*$. 
\begin{defi} \label{defmu}
 Let $F:H\rightarrow G$ be an operator. The \textit{quasi-optimality} constant
 $\mu=\mu(F)>0$ is the smallest number $\mu$, such that 
 \begin{equation}\label{mu}
  \|f-Ff\|\leq \mu \|f-P_G f\|,\quad \textnormal{for all } f\in H,
 \end{equation}
where $P_G:H \rightarrow G$ is the orthogonal projection onto $G$.
If there does not exist a $\mu\in \mathbb{R}$ such that 
\textnormal{(\ref{mu})} is fulfilled, we set $\mu=\infty$. 
\end{defi}
We note that $P_G f$ is the element of $G$ closest to $f$. Thus the \textit{quasi-optimality} constant $\mu(F)$
is a measure of how well $F$ performs in comparison to $P_G$. 

In order to measure stability of the reconstruction, we define the quantity $\eta(F)$ as the operator norm of $Q_{|\R(\U^*)}$. 
\begin{defi} \label{defnu}
 Let $F:H\rightarrow G$ be an operator such that, for each $f\in H$, $Ff$ 
 depends only on the measurements $\U^*f$, i.e., $F=Q\mathcal{U}^*$ for 
 some operator $Q:l^2(\mathbb{N})\rightarrow G$. We define $\eta=\eta(F)$ by 
\begin{equation}
\eta(F):=\underset{\U^*f\neq 0}{\textnormal{sup}}~\frac{\|Q\U^*f\|}{\|\U^*f\|}= \|Q_{|\R(\U^*)}\|.
\end{equation}
\end{defi}
If $\eta(F)$ is small, we call $F$ a well-conditioned mapping, and otherwise ill-conditioned.

In section \ref{etamu} we show that for both, the oblique projection $\mathcal{P}_{G P_U(G)^\perp}$ and for the 
oblique projection $\mathcal{P}_{G \s(G)^\perp}$, the mapping introduced in \cite{1,adcock,adhapo12}, 
$\eta(F) = \|Q\|$ holds. In this case the following lemma applies.

\begin{lem}\label{errestimate}
 Let $F:H\rightarrow G$ be an operator that can be decomposed into the form $F=Q\mathcal{U}^*$ for 
 some operator $Q:l^2(\mathbb{N})\rightarrow G$. Let $f\in H$ and $c\in l^2(\mathbb{N})$. If $\eta(F) = \|Q\|$, then 
 \begin{equation}\label{gerrestimate}
\|P_G f-Q(\mathcal{U}^*f+c)\|\leq \|f-P_G f\|\sqrt{\mu^2-1}+\|c\|\eta. 
 \end{equation}
\begin{proof}
Using the Pythagorean theorem, Definition \ref{defmu} and Definition \ref{defnu} we obtain
 \begin{align*}
 \|P_G f-Q(\mathcal{U}^*f+c)\|&\leq \|P_G f -F f\|+\|Q c\|\\
 &= \sqrt{\|f-Ff\|^2-\|f-P_Gf\|^2}+\eta\|c\| \\
 &\leq \sqrt{\mu^2\|f-P_G f\|^2-\|f-P_Gf\|^2}+\|c\|\eta \\
 &\leq \|f-P_G f\|\sqrt{\mu^2-1}+\|c\|\eta.
 \end{align*}
\end{proof}
\end{lem}

Equation \ttt{gerrestimate} bounds the distance between $P_Gf$ and $F(f)$,
where $F(f)$ is calculated from the perturbed measurements $\mathcal{U}^*f+c$. We can obtain a good error estimate if $f$ is not close to the reconstruction space, provided that $\|c\|$ is small and $\mu$ is close to one.
\begin{defi}
Let $S$ and $W$ be closed subspaces of a Hilbert space $H$ and let $P_W:H\rightarrow W$ be the orthogonal projection onto 
$W$. We define the subspace angle $\varphi = \varphi_{SW}\in [0,\frac{\pi}{2}]$ between $S$ and $W$ by 
\begin{equation}\label{angle}
\cos(\varphi_{SW}) = \underset{\underset{\|s\|=1}{s\in S}}{\inf} \|P_Ws\|.
\end{equation}
\end{defi}

Let $\mathcal{S} = \U \U^*$ denote the frame operator of the frame $\{u_j\}_{j\in \mathbb{N}}$, the frame for the sampling space. In \cite{adhapo12} it is shown for finite dimensional $G$, that if $\cos(\varphi_{GU})>0$, then the oblique projection
with range $G$ and nullspace $\mathcal{S}(G)^\perp$, denoted by $\PP_{G~ \mathcal{S}(G)^\perp}$, exists, and can be written in the form (see \cite[section 4.4]{adhapo12})
\begin{equation} \label{11111}
\PP_{G ~\mathcal{S}(G)^\perp} = \G(\U^*\G)^\dagger \U^*.
\end{equation}
Therefore if 
 $Q_1:= G(\U^*\G)^\dagger $, 
then the oblique projection factors as 
\begin{equation}\nonumber
\PP_{G~ \mathcal{S}(G)^\perp} = Q_1\U^*.
\end{equation}

This shows that the oblique projection $\PP_{G ~\mathcal{S}(G)^\perp}f$ can be calculated from the measurements
$\{\langle f,u_j\rangle\}_{j\in \mathbb{N}}$.

Formula \ttt{11111} is equivalent to
\begin{equation} \nonumber
\PP_{G ~\mathcal{S}(G)^\perp}f = \G\hat{c},
\end{equation}
where $\hat{c}$ is the minimal norm element of the set
 \begin{equation} \label{leastsquare1}
 \underset{c}{\textnormal{arg~min}}\|
  \mathcal{U}^*f-\mathcal{U}^*\mathcal{G}c\|.
 \end{equation}

The following theorem can be found in \cite[Theorem 6.2.]{adhapo12}.
\begin{thm} \label{best}
Let $\{u_j\}_{j=1,\dots, m}$ and $\{g_k\}_{k=1,\dots,n}$ be finite sequences in $H$, and $\cos(\varphi_{GU})>0$.
Let $F:H\rightarrow G$ be a mapping that can be decomposed into $F=Q_2\mathcal{U}^*$ for 
some mapping $Q_2:\R(\mathcal{U}^*)\rightarrow G$. 
If $F(f) = f$ for all $f\in G$, then 
\begin{equation}\nonumber
\eta(F)\geq \eta(\PP_{G ~\mathcal{S}(G)^\perp}).
\end{equation}
\end{thm}
If the quasi-optimality constant $\mu(F)<\infty$, then $F(f) = f$ for all $f\in G$. In this case
Theorem \ref{best} states that $\PP_{G~ \mathcal{S}(G)^\perp}$ has the smallest possible $\eta(F)$ among all 
$F=Q \U^*$.

The main theorems of this paper are Theorems \ref{wichtig} and  \ref{smallestmu}.

\begin{thm} \label{wichtig}
 If $\cos(\varphi_{GU})>0$, then $\mathcal{P}_{GP_U(G)^\perp}$, the oblique projection with range $G$ and kernel
 $P_U(G)$ exists and
 \begin{equation}\nonumber
\mathcal{G}\left((\mathcal{U}^*\mathcal{U})^{\frac{\dagger}{2}}\mathcal{U}^*\mathcal{G}\right )^\dagger
  (\mathcal{U}^*\mathcal{U})^{\frac{\dagger}{2}} \mathcal{U}^* = \mathcal{P}_{GP_U(G)^\perp} .
 \end{equation}
 Equivalently, 
 \begin{equation}\nonumber
 \mathcal{P}_{GP_U(G)^\perp}f=\mathcal{G}\hat{c},
 \end{equation}
 where $\hat{c}$ is the minimal norm element of the set
 \begin{equation} \label{leastsquare}
 \underset{c}{\textnormal{arg~min}}~\|(\mathcal{U}^*\mathcal{U})^{\frac{\dagger}{2}}
  \mathcal{U}^*f-(\mathcal{U}^*\mathcal{U})^{\frac{\dagger}{2}}\mathcal{U}^*\mathcal{G}c\|.
 \end{equation}
\end{thm}
 
\begin{thm}\label{muframeindependent}
If $\cos(\varphi_{GU})>0$, then 
\begin{equation} \label{4.4}
\|\mathcal{P}_{G P_U(G)^\perp}\|= \frac{1}{\cos(\varphi_{GU})}
\end{equation}
and 
\begin{equation} \label{4.5}
\|f-P_G f\|\leq \|f-\mathcal{P}_{G P_U(G)^\perp}f\|\leq \frac{1}{\cos(\varphi_{GU})}\|f-P_G f\|.
\end{equation}
The bound in \ttt{4.5} is sharp.
\end{thm}
 
\begin{thm} \label{smallestmu}
 Let $\cos(\varphi_{GU})>0$. If the mapping $F=Q\U^*:H\rightarrow H$, satisfies
 \begin{equation}\label{jaja}
  \|f-Ff\|\leq\alpha\|f-P_Gf\|, \quad \text{ for all }f\in H,
 \end{equation}
 for some $\alpha>0$, then $\alpha\geq \mu(\mathcal{P}_{GP_U(G)^\perp})$.
\end{thm}

Theorem \ref{wichtig} shows how the oblique projection $\mathcal{P}_{GP_U(G)^\perp}f$ can be calculated from 
the measurements $\{\langle f,u_j\rangle \}_{j\in \mathbb{N}}$. Specifically, setting 
\begin{equation} \label{Q}
 Q:= \mathcal{G}\left((\mathcal{U}^*\mathcal{U})^{\frac{\dagger}{2}}\mathcal{U}^*\mathcal{G}\right )^\dagger
  (\mathcal{U}^*\mathcal{U})^{\frac{\dagger}{2}},
\end{equation}
this projection is given by 
\begin{equation}\nonumber
 \mathcal{P}_{G P_U(G)^\perp}=Q\U^*.
\end{equation}
Theorem \ref{muframeindependent} shows that
the quasi optimality constant of this projection is $\frac{1}{\cos(\varphi_{GU})}$,
and Theorem \ref{smallestmu} states that
this is smallest possible quasi-optimality constant. 

A key property of the mapping $\mathcal{P}_{GP_U(G)^\perp}$ is that $\mu(\mathcal{P}_{GP_U(G)^\perp})$ and $\eta(\mathcal{P}_{GP_U(G)^\perp})$ can be calculated. In section \ref{etamu} we state explicit formulas for them.

\section{Existence of the oblique projection}
The following lemma follows from \cite[Thm. 2.1]{2000} and \cite[(2.2)]{2000}.
\begin{lem} \label{l}
 If $S$ and $W$ are closed subspaces of a Hilbert space $H$, then $\cos(\varphi_{SW^\perp})>0$ if and only if
 $S\cap W=\{0\}$ and $S\oplus W$ is closed in $H$.
\end{lem} 

The proof of the following lemma can be found in \cite[Theorem 1]{ex}.
\begin{lem} \label{exist}
If $S \cap W = \{0\}$ and $H_1:=S\oplus W$ is a closed subspace of $H$, then the
oblique projection $\mathcal{P}_{SW}:H_1\rightarrow S$ with range $S$ and kernel $W$ is well defined and bounded.
\end{lem}

The following theorem can be found in \cite[Corollary 3.5]{adhapo12}.
\begin{thm} \label{coroll}
 Let $S$ and $W$ be closed subspaces of $H$ with $\cos(\varphi_{SW^\perp})>0$ and let $H_1:=S\oplus W$. If
 $\mathcal{P}_{SW}:H_1\rightarrow S$ is the oblique projection with range $S$ and kernel $W$, then
 \begin{equation} \label{aaa}
  \|\mathcal{P}_{SW}\|=\frac{1}{\cos(\varphi_{SW^\perp})}
 \end{equation}
 and
\begin{equation} \label{bbb}
 \|f-P_Sf\|\leq\|f-\mathcal{P}_{SW}f\|\leq \frac{1}{\cos(\varphi_{SW^\perp})}\|f-P_Sf\|,
\end{equation}
for all $f\in H_1$.
The upper bound in \textnormal{(\ref{bbb})} is sharp.
\end{thm}

We make use of the following well known lemma.
\begin{lem} \label{closedrange}
Let $L$ and $H$ be Hilbert spaces, and let $\U:H\rightarrow L$ be a bounded operator. If there exists an $A>0$, such that
\begin{equation} \label{unt}
 A\|c\|\leq\|\U c\| \textnormal{ for all }c\in \N(\U)^\perp, 
\end{equation}
then the operator $\U$ has a closed range.
\end{lem}
In the following, we use the notation
\begin{equation}\nonumber
 P_U(G):=\{P_Ug:g\in G\}.
 \end{equation}
 
\begin{lem}\label{lemma}
Let $G$ and $U$ be closed subspaces of a Hilbert space $H$. If $\cos(\varphi_{GU})>0$, then the subspace
$P_U(G)$ is closed. Furthermore
 \begin{equation}\label{cosvergleich}
 \cos(\varphi_{G P_U(G)})=\cos(\varphi_{GU}).
 \end{equation}
 \begin{proof}
 
 From \textnormal{(\ref{angle})}, it follows that
 \begin{equation}\nonumber
 \|g\|\cos(\varphi_{GU})\leq \|P_Ug\|\quad \textnormal{for all }g\in G.
 \end{equation}
 The closedness of the subspace $P_U(G)$ follows from the fact that $\N(P_U)^\perp = (U^\perp)^\perp = \overline{U} = U$,
using Lemma \ref{closedrange}.
 The second statement follows from
\begin{align*}
&\cos(\varphi_{G P_U(G)})=\underset{\underset{\|u\|=1}{u\in G}}{\inf}~ \underset{\underset{\|v\|=1}{v\in P_U(G)}}{\sup} 
\langle u,v\rangle=\underset{\underset{\|u\|=1}{u\in G}}{\inf} ~\underset{\underset{\|v\|=1}{v\in P_U(G)}}{\sup} 
\langle P_Uu,v\rangle=\\
&\underset{\underset{\|u\|=1}{u\in G}}{\inf}\|P_Uu\|=\cos(\varphi_{GU}),
\end{align*}
using \textnormal{(\ref{angle})} for the last equality.
\end{proof}
\end{lem}

\begin{thm} \label{lll}
If $\cos(\varphi_{GU})>0$, then $H=G \oplus P_U(G)^\perp$, and the oblique projection 
$\mathcal{P}_{G P_U(G)^\perp}:H\rightarrow G$ is well defined and bounded.
\begin{proof}
By assumption $\cos(\varphi_{GU})>0$ and thus by Lemma \ref{lemma} 
\begin{equation}\nonumber
\cos(\varphi_{G P_U(G)})=\cos(\varphi_{GU})>0.
\end{equation}
By Lemma \ref{l} and Lemma \ref{exist}
the oblique projection
$\mathcal{P}_{G P_U(G)^\perp}$ is well defined and bounded as a mapping
from $G \oplus P_U(G)^\perp$ onto $G$. We prove that 
\begin{equation}\nonumber
G\oplus P_U(G)^\perp = H.
\end{equation}
Lemma \ref{l} implies that $G\oplus P_U(G)^\perp$ is closed. Consequently, it is sufficient to show that 
${(G\oplus P_U(G)^\perp)}^\perp = \{0\}$.
By assumption $\cos(\varphi_{GU})>0$ and thus by Lemma \ref{lemma} $P_U(G)$ is closed, and
\begin{equation}\nonumber
 {(G\oplus P_U(G)^\perp)}^\perp = G^\perp\cap \overline{P_U(G)}=G^\perp \cap P_U(G).
\end{equation}
Let $h\in G^\perp \cap P_U(G)$. Using that $h=P_U g$ for some $g\in G$ and that $h\in G^\perp$, we conclude that
for all $s\in G$
\begin{equation} \label{lax}
0 = \langle h,s\rangle = \langle g,P_U s\rangle.
\end{equation}
If $g\neq 0$, from \ttt{lax}, it follows that 
\begin{equation}\nonumber
 \underset{\underset{\|g\|=1}{g\in G}}{\inf}~ \underset{\underset{\|P_U s\|=1}{s\in G}}{\sup}~ \left\langle g,P_U s\right\rangle = 0.
\end{equation}
This is a contradiction to $\cos(\varphi_{G P_U(G)})>0$. Consequently $g = 0$ and  $h = P_U g = 0$.
\end{proof}
\end{thm}
\begin{lem} \label{schnitt}
Assume that $\{u_j\}_{j\in \mathbb{N}}$ is a frame sequence in $H$ and $\{g_k\}_{k\in \mathbb{N}}$ is a Riesz sequence in $H$. Then
$G\cap U^\perp = \{0\}$ if and only if the operator $\U^*\G$ is injective.
\begin{proof}
Let $G\cap U^\perp = \{0\}$. This implies that for every $g\in G$ with $\U^*g = 0$, it follows that $g = 0$.
For every $g \in G$ there exists a $c\in l^2(\mathbb{N})$ such that $g = \G c$, and consequently for every 
$c \in l^2(\mathbb{N})$ with $\U^* \G c = 0$, it follows that $\G c = 0$. Since $\{g_k\}_{k\in \mathbb{N}}$ is a Riesz sequence,
$\G c = 0$ if and only if $c= 0$, which shows that the operator $\U^*\G$ is injective.

Let the operator $\U^*\G$ be injective. Using that $\{g_k\}_{k\in \mathbb{N}}$ is a Riesz sequence, we conclude that
for every $c \in l^2(\mathbb{N})$ with $\U^*\G c = 0$, it follows that $\G c = 0$. Since 
$\{g_k\}_{k\in \mathbb{N}}$ is a frame for $G$, this implies that for every 
$g \in G$ with $\U^*g = \{0\}$, it follows that $g = 0$. Consequently $G\cap U^\perp = \{0\}$.
\end{proof}
\end{lem}

Lemma \ref{schnitt} implies, that for finite sequences $\{u_j\}_{j\in J}$ and $\{g_k\}_{k\in K}$ with $g_k$, $k\in K$, linearly 
independent, $\cos(\varphi_{GU})>0$ if and only if $\U^*\G$ is injective.

\section{Frames  and  the Pseudoinverse}
We make use of the following version of the spectral theorem.
\begin{thm}\label{spectral}
Let $A$ be a bounded operator on $H$. If $A$ is normal, there exists a measure space $(X,\Sigma,\mu)$
and a real-valued essentially bounded measurable function $f$ on $X$ and a unitary operator
$U:H\rightarrow L^2_\mu(X)$ such that 
\begin{equation} \label{gspectral}
U^*TU = A,
\end{equation}
where 
\begin{equation}\nonumber
[T g](x) = f(x) g(x) 
\end{equation}
\begin{equation}\nonumber
\text{ and }\|T\| = \|f\|_\infty
\end{equation}
\end{thm}
For an essentially bounded function $f\in L^2_\mu(X)$ we use the notation $M_f$ for the multiplication operator  
\begin{equation}\nonumber
M_fg(x):= f(x)g(x)\text{ for }g\in L^2_\mu(X).
\end{equation}
The following theorem can be found in \cite[Theorem 2.1]{cr}
\begin{thm} \label{crange}
 Let $ f \in L^2_\mu(X)$ be essentially bounded. Then $M_f$ has a closed range if and only if $f$ is bounded
 away from zero on $X\backslash \{x\in X:f(x) = 0\}$.
\end{thm}

We need the definition of the pseudoinverse in a Hilbert space.
\begin{lem}\label{pseudoinverse}
Let $H$ and $L$ be Hilbert spaces. If $\U:L\rightarrow H$ is a bounded operator with a closed range $\R(\U)$, then
there exists a unique bounded operator $\U^\dagger:H\rightarrow L$ such that 
\begin{align}
 &\N(\U^\dagger) = \R(\U)^\perp,\label{pi1} \\ 
 &\R(\U^\dagger) = \N(\U)^\perp, \text{ and }\label{pi2}\\ 
 &\U \U^\dagger x =x,~x\in \R(\U).\label{pi3}
 \end{align}
\end{lem}
We call the operator $A^\dagger$ the pseudoinverse of $A$.
The proof of the following lemma is straightforward and thus skipped. 
\begin{lem} \label{straight}
Let $f \in L^2_\mu(X)$ be an essentially bounded function. If $f$ is bounded away from zero 
on $X\backslash \{x\in X:f(x) = 0\}$, then the pseudoinverse of the multiplication operator $M_f$ is a bounded
function and given by
\begin{equation}\nonumber
 M_f^\dagger g(x) = 
 \begin{cases}
 \frac{g(x)}{f(x)} \quad &\text{for } x\in  X\backslash \{x\in X:f(x) = 0\},\\
 0 &\text{otherwise. }
 \end{cases}
\end{equation}
\end{lem}

\begin{lem} \label{nr}
 If $\U:L\rightarrow H$ is an operator with a closed range, then the following identities hold,
 \begin{align*}
 &\R(\U\U^*) = \R(\U), &\N(\U\U^*) = \N(\U^*),\\
 &\R(\U^*\U) = \R(\U^*), &\N(\U^*\U)=\N(\U).
 \end{align*}
\begin{proof}
 Obviously $\R(\U\U^*) \subset \R(\U)$. Let $f\in \R(\U)$, $f= \U d$ for some $d\in L$.
 Since  $\R(\U^*)$ is closed ($\R(\U)$ is closed if and onlay if $\R(\U)$ is closed), $d$ can be decomposed into 
 $d = d_{\R(\U^*)}+d_{\R(\U^*)^\perp}$, where 
 $d_{\R(\U^*)}=P_{\R(\U^*)}c \in \R(\U^*)=\overline{\R(\U^*)} = \N(\U)^\perp$ and 
 $d_{\R(\U^*)^\perp}= d-P_{\R(\U^*)}d\in \R(\U^*)^\perp=\N(\U)$. Consequently $f = \U d_{\R(\U^*)}$ and 
 $f\in \R(\U \U^*)$.
 The proofs of the other statements are similar.
\end{proof}
\end{lem}

We observe that \ttt{frameequation} can be written in the form 
\begin{equation}\nonumber
 A \leq \frac{\langle \s f,f\rangle}{\langle f,f\rangle} \leq B, \quad \text{for every }f\in \R(\U). 
\end{equation}
Since by Lemma \ref{nr} $\R(\U) = \N(\U^*)^\perp = \N(\U \U^*)^\perp =
\N(\s)^\perp$, this ensures that, except of zero, the spectrum of the operator $\s$ is bounded away from zero. 
Using Theorem \ref{spectral}, Theorem \ref{crange} and Lemma \ref{straight} we obtain that the pseudoinverse of
the operator $\s^\dagger$ exists and is a bounded operator on $H$.
Since $\s\geq 0$ also $\s^\dagger\geq 0$. For every positive operator there exists a unique positive
square root. Therefore we can define the operator
\begin{equation}\nonumber
 \s^{\frac{\dagger}{2}}:= {(\s^\dagger)}^\frac{1}{2}.
\end{equation}

The following lemma can be found in \cite[Lemma 5.4.5]{ch08}
\begin{lem}\label{545}
Let $\{u_j\}_{j\in \mathbb{N}}$ be a sequence in H. The sequence $\{u_j\}_{j\in \mathbb{N}}$
is a frame sequence with frame bounds $A$ and $B$ if and only if the synthesis operator $\U$ is well defined on $l^2(\mathbb{N})$ and
\begin{equation} \label{e545}
A \|c\|^2\leq \|\U c\|^2\leq B \|c\|^2\quad \text{for all } 
c\in {\mathcal{N}(\U)}^\perp
\end{equation}
\end{lem}

We observe, that \ttt{e545} can be written in the form
\begin{equation}\nonumber
 A \leq \frac{\langle \U^*\U c ,c\rangle}{\langle c,c\rangle} \leq B\quad  \text{for every }c\in {\mathcal{N}(\U)}^\perp.
\end{equation}
Since by Lemma \ref{nr} ${\mathcal{N}(\U)}^\perp = {\mathcal{N}(\U^*\U)}^\perp$, this ensures that except of zero, the spectrum of the operator $\U^*\U$ is bounded away from zero. 
Using Theorem \ref{spectral}, Theorem \ref{crange} and Lemma \ref{straight} we obtain that the pseudoinverse of
the operator $(\U^*\U)^\dagger$ exists and is a bounded operator on $l^2(\mathbb{N})$.
Since $\U^*\U\geq 0$ also $(\U^*\U)^\dagger\geq 0$. For every positive operator there exists a unique positive
square root. Therefore we can define the operator
\begin{equation}\nonumber
 (\U^*\U)^{\frac{\dagger}{2}}:= \Big((\U^*\U)^\dagger \Big) ^\frac{1}{2}.
\end{equation}

The following theorem is a slightly modified version of \cite[Theorem 5.3.4]{ch08} and thus omitted.

\begin{thm} \label{tight}
Let $\{u_j\}_{j\in\mathbb{N}}$ be a frame sequence in $H$. If $\mathcal{S}$ is the corresponding frame operator, then 
\begin{equation}\nonumber
 \{\mathcal{S}^{\frac{\dagger}{2}} u_j\}_{j\in\mathbb{N}}
\end{equation}
forms a tight frame for $U$ with frame bound equal to $1$. Let $M = \s^{\frac{\dagger}{2}}\U$ denote the 
synthesis operator of the tight frame sequence $\{\mathcal{S}^{\frac{\dagger}{2}} u_j\}_{j\in\mathbb{N}}$.
Then 
\begin{equation} \label{p3}
 P_U = MM^* = \mathcal{S}^{\frac{\dagger}{2}}\mathcal{S}\mathcal{S}^{\frac{\dagger}{2}}.
\end{equation}
\end{thm}

\begin{lem}\label{tight1}
Let $\{u_j\}_{j\in\mathbb{N}}$ be a frame sequence in $H$.
The operator
 $(\U^*\U)^{\frac{\dagger}{2}}\U^*$
is the analysis operator of the frame sequence $\{\s^{\frac{\dagger}{2}} u_j\}_{j\in\mathbb{N}}$.
Equivalently, 
\begin{equation} \label{hola}
 (\U^*\U)^{\frac{\dagger}{2}}\U^* = \U^*(\U\U^*)^{\frac{\dagger}{2}}.
 \end{equation}
\begin{proof}
 Obviously for $k\in \mathbb{N}$
 \begin{equation}\nonumber
 (\U^*\U)^k\U^* = \U^*(\U \U^*)^k. 
 \end{equation}
Therefore 
\begin{equation}\nonumber
 \gamma(\U^*\U)\U^* = \U^*\gamma(\U \U^*) 
\end{equation}
 for every polynomial $\gamma$. Taking limits, it follows that 
 \begin{equation}\nonumber
 f(\U^*\U)\U^* = \U^*f(\U \U^*) 
\end{equation}
for every continous function $f$, in particular for 
\begin{equation}\nonumber
 f(A) = A^\frac{\dagger}{2}.
\end{equation}
\end{proof}
\end{lem}

\begin{lem} \label{3021}
  If $\cos(\varphi_{GU})>0$, then 
  \begin{equation} \label{null2}
  \mathcal{N}(\G) = \N((\U^*\U)^{\frac{\dagger}{2}}\U^* \G),
  \end{equation}
  and
\begin{equation} \label{cond}
\sqrt{C}\cos(\varphi_{GU})\|c\|\leq \|(\U^*\U)^{\frac{\dagger}{2}}\U^* \G c\|\leq \sqrt{D}\|c\| \quad 
\text{for all }c\in {\mathcal{N}(\G)}^\perp. 
\end{equation}
\begin{proof}
We define 
\begin{equation} \label{defL}
L:=\U(\U^*\U)^{\frac{\dagger}{2}}.
\end{equation}
 By Theorem  \ref{tight1}, the operator  $L^*$ is the analysis operator 
 of the tight frame sequence $\{\mathcal{S}^{\frac{\dagger}{2}} u_j\}_{j\in\mathbb{N}}$ with frame bound equal to one.
 Therefore
\begin{equation}\label{io} 
\|L^*f\|=\|f\| \quad \text{for all }f\in U.
\end{equation}

Clearly $P_U \U = \U$ and consequently
\begin{equation} \label{gup}
\U^* = \U^*P_U.
\end{equation}
Using \ttt{gup}, \textnormal{(\ref{io})} and $\textnormal{(\ref{e545})}$, we obtain 
\begin{equation} \label{upper}
\|L^* \G c\|= \|L^* P_U(\G c)\|=\|P_U(\G c)\|\leq \|\G c\|\leq \sqrt{D}\|c\|.
\end{equation}
Using the definition of $\cos(\varphi_{GU})$ and $\ttt{e545}$, it follows that for $c\in {\mathcal{N}(\G)}^\perp$
\begin{equation}\label{lower}
\|P_U(\G c)\|\geq \cos(\varphi_{GU}) \|\G c\|\geq  \cos(\varphi_{GU})\sqrt{C} \|c\|,
\end{equation}
which proves \textnormal{(\ref{cond})}.

Trivially $\N(\G)\subset  \N(L^* \G)$. 
Let $c\in \mathcal{N}(L^* \G)$.
We decompose $c$ into $c = c_1+c_2$, with $c_1:=P_{\N(G)}c\in \N(G)$ and $c_2:=c-P_{\N(G)}c\in \N(G)^\perp$.
Using \textnormal{(\ref{cond})}, we obtain that 
\begin{equation}\nonumber 
0 = \|L^* \G c\|=\|L^* \G c_2\|\geq \sqrt{C}\cos(\varphi_{GU})\|c_2\|,
\end{equation}
which implies that $c_2=0$, and consequently $c\in \N(G)$.
\end{proof}
\end{lem}

The following lemma is a part of \cite[Lemma 2.5.2]{ch08}.
\begin{lem}\label{252}
Let $\U:L\rightarrow H$ be a bounded operator. If $\U$ has a closed range $\R(\U)$, then the following holds:
\begin{enumerate}
\item The orthogonal projection of $H$ onto $\R(\U)$ is given by $\U\U^\dagger$. \label{pl1}
\item The orthogonal projection of $L$ onto $\R(\U^\dagger)$ is given by $\U^\dagger \U$.\label{pl2}
\item  The operator $\U^*$ has closed range, and $(\U^*)^\dagger = (U^\dagger)^*$.
\end{enumerate}
\end{lem}

The following Lemma can be found in \cite[Corollary 1.]{pseudoinv}
\begin{lem}\label{cns}
If $\U:L\rightarrow H$ is a bounded operator with a closed range, then 
\begin{equation}\label{nullconj}
\N(\U^\dagger) = \N(\U^*)  
\end{equation}
and
\begin{equation}\label{rangeconj}
\R(\U^\dagger) = \R(\U^*).  
\end{equation}
\begin{proof}
 By \ttt{pi1}, 
 \begin{equation} \nonumber
 \N(\U^\dagger)=\R(\U)^\perp=\N(\U^*).
 \end{equation}
 Since $\R(\U)$ is closed, if and only if $\R(\U^*)$ is closed, by \ttt{pi2},
 \begin{equation} \nonumber
 \R(\U^\dagger)=\N(\U)^\perp=\overline{\R(\U^*)}=\R(\U^*).
 \end{equation}
\end{proof}
\end{lem}

\section{Proofs}

\textbf{Proof of  Theorem \ref{wichtig}.}
We set $R:=\mathcal{G}\left((\mathcal{U}^*\mathcal{U})^{\frac{\dagger}{2}}\mathcal{U}^*\mathcal{G}\right )^\dagger
(\mathcal{U}^*\mathcal{U})^{\frac{\dagger}{2}} \U^*$. We show 
\begin{enumerate}
 \item $R$ is well defined,
 \item $\R(R) = G$,
 \item $\N(R) = P_U(G)^\perp$ and
 \item $R^2 = R$.
\end{enumerate}

(1.) The lower bound of \textnormal{(\ref{cond})} ensures that the operator $(\U^*\U)^{\frac{\dagger}{2}}\U^* \G$
has a closed range. The upper bound of \ttt{cond} shows that the operator
$(\U^*\U)^{\frac{\dagger}{2}}\U^* \G$ is bounded and by Lemma \ref{pseudoinverse}
this proves the existence of 
the operator $\left((\mathcal{U}^*\mathcal{U})^{\frac{\dagger}{2}}\mathcal{U}^*\mathcal{G}\right )^\dagger$.
Therefore $R$ is a well defined operator from $H$ to $G$.\\
(2.) Clearly $\R(R)\subset G$. The lower bound in Lemma \ref{3021} implies that $\R((\mathcal{U}^*\mathcal{U})^{\frac{\dagger}{2}}\mathcal{U}^*\mathcal{G})$ is closed.
From Lemma \ref{252} (\ref{pl2}.) and Lemma \ref{cns}, it follows that 
\begin{align*}
R\G &= \G \left((\mathcal{U}^*\mathcal{U})^{\frac{\dagger}{2}}\mathcal{U}^*\mathcal{G}\right )^\dagger
(\mathcal{U}^*\mathcal{U})^{\frac{\dagger}{2}} \mathcal{U}^*\mathcal{G}
= \G P_{\overline{\R\left(((\mathcal{U}^*\mathcal{U})^{\frac{\dagger}{2}}\mathcal{U}^*\mathcal{G})^*\right)}}\\
& = \G P_{\mathcal{N}\left((\mathcal{U}^*\mathcal{U})^{\frac{\dagger}{2}}\mathcal{U}^*\mathcal{G}\right)^\perp} = 
\G P_{\N (\G)^\perp}.
\end{align*}
(3.) Let $L$ be defined by \ttt{defL}. From Theorem \ref{tight} and Lemma \ref{tight1}, we deduce that
$P_U = L L^*= \U(\U^*\U)^\dagger \U^*$. Since $\{g_k\}_{k\in \mathbb{N}}$ is a frame for $G$ and $\G$ is the corresponding synthesis operator, it holds $\R(\G) = G$ and
consequently
\begin{align}\label{abc}
 P_U(G)^\perp =& \R \left (\U(\U^*\U)^\dagger \U^*\G\right )^\perp =\N\left ([\U(\U^*\U)^\dagger \U^*\G]^*\right)\\
 =& \N(\G^*\U(\U^*\U)^\dagger \U^*).
\end{align}
From \ttt{nullconj}, it follows that
\begin{equation}\label{defg}
\mathcal{N}\left(((\mathcal{U}^*\mathcal{U})^{\frac{\dagger}{2}}\mathcal{U}^*\mathcal{G} )^\dagger
  (\mathcal{U}^*\mathcal{U})^{\frac{\dagger}{2}} \mathcal{U}^*\right)=
  \N(\G^*\U (\U^*\U)^\dagger \U^*).
\end{equation}
Combining \ttt{abc} and \ttt{defg}, we obtain
\begin{equation}\nonumber
P_U(G)^\perp\subset \N \left(\mathcal{G}\left((\mathcal{U}^*\mathcal{U})^{\frac{\dagger}{2}}\mathcal{U}^*\mathcal{G}\right )^\dagger
(\mathcal{U}^*\mathcal{U})^{\frac{\dagger}{2}} \U^*\right)= \N(R).
\end{equation}
In Theorem \ref{lll}, it is shown that $G\oplus P_U(G)^\perp =H$. 
Let $f\in \N(R)$.
We decompose $f$ into $f = f_1 + f_2$, where $f_1\in G$ and $f_2\in P_U(G)^\perp$. Since $f\in \N(R)$ and $P_U(G)^\perp\subset \N(R)$,
\begin{equation}\nonumber
 0 = Rf = R(f_1+f_2) = R f_1 = f_1,
\end{equation}
and consequently $f = f_2 \in P_U(G)^\perp$.

(4.)
Using Lemma \ref{252} \textnormal{(\ref{pl2}.)}, we obtain 
\begin{align*}
R^2  &=\mathcal{G}\left((\mathcal{U}^*\mathcal{U})^{\frac{\dagger}{2}}\mathcal{U}^*\mathcal{G}\right )^\dagger
(\mathcal{U}^*\mathcal{U})^{\frac{\dagger}{2}} \U^*\mathcal{G}\left((\mathcal{U}^*\mathcal{U})^{\frac{\dagger}{2}}\mathcal{U}^*\mathcal{G}\right )^\dagger
(\mathcal{U}^*\mathcal{U})^{\frac{\dagger}{2}} \U^*\label{idempotent}\\
&= \mathcal{G} P_{\R\left(((\mathcal{U}^*\mathcal{U})^{\frac{\dagger}{2}}\mathcal{U}^*\mathcal{G})^\dagger \right)}                  \left((\mathcal{U}^*\mathcal{U})^{\frac{\dagger}{2}}\mathcal{U}^*\mathcal{G}\right )^\dagger
(\mathcal{U}^*\mathcal{U})^{\frac{\dagger}{2}} \mathcal{U}^*.\nonumber                   \\ 
&= \mathcal{G}\left((\mathcal{U}^*\mathcal{U})^{\frac{\dagger}{2}}\mathcal{U}^*\mathcal{G}\right )^\dagger
(\mathcal{U}^*\mathcal{U})^{\frac{\dagger}{2}} \mathcal{U}^* = R.\nonumber
\end{align*}

\textbf{Proof of Theorem \ref{muframeindependent}.} The statement
follows by combining 
\ttt{cosvergleich} and 
 Theorem \ref{coroll}. 

\textbf{Proof of  Theorem \ref{smallestmu}.}
Let $g\in G$. Every element  $f\in g+U^\perp$ has the same value $\U^*g$, and thus 
 \begin{equation}\label{jaja1}
 Ff=Q\U^*g \quad \text{for all }f\in g+U^\perp.
 \end{equation}
 From \ttt{jaja}, it follows that 
 \begin{equation}\label{jaja2}
 Fh=h \quad \text{ for all }h\in G, 
 \end{equation}
 otherwise $\mu = \infty$, and consequently $Q\U^*g = g$.
 From \ttt{jaja1}, we deduce that 
 \begin{equation}\nonumber
  F(g+u^\perp) = g \quad \text{for all }g\in G \text{ and } u\in U^\perp.
 \end{equation}
This means that $F_{|G\oplus U^\perp}=\mathcal{P}_{GU^\perp}$. From \ttt{bbb}, it follows that
\begin{equation}\nonumber
\mu(F_{|G\oplus U^\perp})=\frac{1}{\cos(\varphi_{GU})},
\end{equation}
which implies that 
$\mu(F)\geq \frac{1}{\cos(\varphi_{GU})}$.
 
\section{An abstract definition of our reconstruction}
The oblique projection $\mathcal{P}_{GP_U(G)^\perp}$ is characterized
as follows. 
Most part of this proof is similar to the proof of \cite[Theorem 4.2]{adhapo12}.
\begin{thm} \label{theorm}
Let $\cos(\varphi_{GU})>0$.
The mapping $\mathcal{P}_{G P_U(G)^\perp}$ is the unique
operator $F$ that satisfies the equations 
\begin{equation} \label{figs}
\langle P_U Ff,g_k\rangle = \langle P_U{f},g_k\rangle, \quad j\in \mathbb{N},~f\in H.
\end{equation}
\begin{proof}
In Theorem \ref{lll} it is shown that $H=G \oplus P_U(G)^\perp$ and that the oblique projection 
$F= \mathcal{P}_{G P_U(G)^\perp}:H\rightarrow G$ is well defined and bounded.
We show next that $\mathcal{P}_{G P_U(G)^\perp}$ satisfies equation 
\textnormal{(\ref{figs})}.
From the self adjointness of $P_U$, and the fact that $\{g_j\}_{j\in \mathbb{N}}$ is a frame sequence for $G$,
it follows that \textnormal{(\ref{figs})} is equivalent to 
\begin{equation} \label{figs1}
\langle Ff,\Phi\rangle=\langle f,\Phi\rangle \quad \text{for all }
\Phi \in P_U(G),~f \in H.
\end{equation}
We have to show that
\begin{equation}\nonumber
\langle \mathcal{P}_{G P_U(G)^\perp}f,\Phi\rangle=\langle f,\Phi\rangle
\quad \text{for all }\Phi \in P_U(G),~ f\in H. 
\end{equation}
Using that $H=G\oplus P_U(G)^\perp $, every $f\in H$ can be decomposed into $f=f_G+f_{P_U(G)^\perp}$, where
$f_G\in G$ and $f_{P_U(G)^\perp}\in P_U(G)^\perp$. Thus
\begin{equation}\nonumber
\langle \mathcal{P}_{G P_U(G)^\perp}(f_G+f_{P_U(G)^\perp}),\Phi
\rangle = \langle f_G,\Phi\rangle = \langle f_G+f_{P_U(G)^\perp},\Phi\rangle = \langle f,\Phi\rangle.
\end{equation}
Next we show the uniqueness. 
We assume that there are two mappings $F_1,F_2:H\rightarrow G$ that satisfy \textnormal{(\ref{figs1})}.
This means for all $f\in H$ and $\Phi\in P_U(G)$
\begin{equation}\label{19}
 \langle F_1f,\Phi\rangle=\langle f,\Phi \rangle= \langle F_2f,\Phi\rangle.
\end{equation}
From \textnormal{(\ref{19})}, it follows that $\R(F_1-F_2)\subset P_U(G)^\perp$.
We know that $\R(F_1)\subset G$ and that $\R(F_2)\subset G$ and thus $\mathcal{R}(F_1-F_2)\subset G\cap 
P_U(G)^\perp$. From Lemma \ref{lemma} in combination with Lemma \ref{l}, it follows that 
$G\cap P_U(G)^\perp=\{0\}$, and consequently $F_1=F_2$.
\end{proof}
\end{thm}

\section{Stability and quasi optimality of $\mathcal{P}_{G P_U(G)^\perp}$}\label{etamu}
In this section we give formulas for the calculation of $\eta(\mathcal{P}_{G P_U(G)^\perp})$ and $\mu(\mathcal{P}_{G P_U(G)^\perp})$.
We also estimate the condition number of the 
operator $(\U^*\U)^{\frac{\dagger}{2}}\U^*\G $ (see \ttt{leastsquare}) in terms of the condition number of $\G$ and $\cos(\varphi_{GU})$.

The following theorem is similar to \cite[Lemma 2.13]{adcock}, but for the convenience  we include a proof. 
\begin{thm}
Let $\{u_j\}_{j=1}^n$ and $\{g_k\}_{k=1}^m$ be finite sequences in $H$.
If the vectors $g_k$, $k=1,\dots m$, are linearly independent, then
\begin{equation}\nonumber
 \cos^2(\varphi_{GU})= \lambda_{\min}((\G^*\G)^{-1}\G^*\U\U^\dagger\G).
 \end{equation}
 \begin{proof}
 By definition
 \begin{equation}\nonumber
 \cos(\varphi_{GU})=\underset{c\neq 0}{\inf}\frac{\|P_U(\G c)\|}{\|\G c\|}.
 \end{equation}
 From $P_U^2=P_U$ and the self adjointness of $P_U$, it follows that
 \begin{equation}\label{rrr}
 \left(\frac{\|P_U(\G c)\|}{\|\G c\|}\right)^2=\frac{\langle\G^* \U \U^\dagger \G c,c\rangle}{\langle\G^*\G c,c\rangle}.
 \end{equation}
 Since the vectors $g_k$, $k=1,\dots m$, are linearly independent $\G^*\G$ is invertible, and so is $(\G^*\G)^{\frac{1}{2}}$. Substituting
 $a:=(\G^*\G)^{\frac{1}{2}} c$ in \ttt{rrr}, we obtain
 \begin{equation}\nonumber
  \left(\frac{\|P_U(\G c)\|}{\|\G c\|}\right)^2=
  \frac{\langle(\G^*\G)^{-\frac{1}{2}}\G^* \U \U^\dagger \G (\G^*\G)^{-\frac{1}{2}}a,a\rangle}{\langle a,a\rangle }.
 \end{equation}
From the self adjointness of $(\G^*\G)^{-\frac{1}{2}}\G^* \U \U^\dagger \G (\G^*\G)^{-\frac{1}{2}}$ it follows that
\begin{equation}\nonumber
 \cos^2(\varphi_{GU})= \lambda_{\min}((\G^*\G)^{-\frac{1}{2}}\G^* \U \U^\dagger \G (\G^*\G)^{-\frac{1}{2}}).
\end{equation}
The two operators $(\G^*\G)^{-\frac{1}{2}}\G^* \U \U^\dagger \G (\G^*\G)^{-\frac{1}{2}}$ and $(\G^*\G)^{-1}\G^* \U \U^\dagger \G$
have the same spectrum, because they are similar. This finishes the proof.
\end{proof}
\end{thm}

\begin{thm}
Set  $Q= \mathcal{G}\left((\mathcal{U}^*\mathcal{U})^{\frac{\dagger}{2}}\mathcal{U}^*\mathcal{G}\right )^\dagger
  (\mathcal{U}^*\mathcal{U})^{\frac{\dagger}{2}}$ and $Q_1:=
  G(\U^*\G)^\dagger$.  Then
\begin{align}
 &\eta(\PP_{G P_U(G)^\perp}) = \|Q_{|\R(\U^*)}\| = \|Q \|,
\text{ and} \label{uujj}\\
&\eta(\PP_{G \s(G)^\perp}) = \|{Q_1}_{|\R(\U^*)}\| = \|Q_1 \|,\label{uujj1}
\end{align}
\begin{proof}
Since by Theorem~\ref{wichtig} $\mathcal{P}_{GP_U(G)^\perp} =  Q
\mathcal{U}^*$, we have  
\begin{align*}
\eta(\PP_{G P_U(G)^\perp})&=\underset{\U^*f\neq 0}{\textnormal{sup}}~\frac{\|Q\U^*f\|}{\|\U^*f\|}
=\|Q_{|\R(\U^*)}\|.
\end{align*}
Since $\R(\U^*)$ is closed, every element in $c\in l^2(\mathbb{N})$ can be decomposed into $c_{\R(\U^*)}+c_{\R(\U^*)^\perp}$ with $c_{\R(\U^*)}: = P_{\R(\U^*)}c\in \R(\U^*)$ and $c_{\R(\U^*)^\perp}:= c-P_{\R(\U^*)}c\in \R(\U^*)^\perp$. If we can prove 
 $\R(\U^*)^\perp\subset \N(Q)$, then  \ttt{uujj} follows. We show that $\R(\U^*)^\perp= \N\left((\U^*\U)^\frac{\dagger}{2}\right)$.
Using Lemma \ref{nr} and \ttt{pi1} we obtain
\begin{equation}\nonumber
 \N\left((\U^*\U)^\frac{\dagger}{2}\right)= \N\left((\U^*\U)^\dagger\right)=\R(\U^*\U)^\perp = \R(\U^*)^\perp.
\end{equation}
The proof of \ttt{uujj1} is similar.
\end{proof}
\end{thm}
For the calculation of the coefficients $\hat{c}$ of the least squares
problem  \ttt{leastsquare}, it is important to know the condition number of the 
operator $(\U^*\U)^{\frac{\dagger}{2}}\U^*\G $. The following
statement gives some hints.  

\begin{thm} \label{5}
If $\cos(\varphi_{GU})>0$, then
\begin{equation} \label{aio}
\cos(\varphi_{GU})\|\G c\|\leq \|(\U^*\U)^{\frac{\dagger}{2}}\U^* \G c\|\leq \|\G c\|\quad \text{for all }c\in l^2(\mathbb{N}),
\end{equation}
and
\begin{equation}\label{aio1}
\cos(\varphi_{GU})\kappa(\G)\leq
 \kappa((\mathcal{U}^*\mathcal{U})^{\frac{\dagger}{2}}\mathcal{U}^*\mathcal{G})
 \leq \frac{1}{\cos(\varphi_{GU})}\kappa(\G).
 \end{equation}
\begin{proof}
Equation \ttt{aio} follows from
\ttt{upper} and \ttt{lower}. 
Equation \ttt{aio1} is a direct consequence of \ttt{aio}.
\end{proof}
\end{thm}

The following Theorem and its proof is similar to \cite[Corollary 4.7]{adhapo12}
\begin{thm} \label{mukond}
Let $\cos(\varphi_{GU})>0$. If $Q$ is defined by \ttt{Q}, then
\begin{equation} \label{21}
\frac{1}{\sqrt{B}}\leq \eta(\mathcal{P}_{G P_U(G)^\perp}) = \|Q\|\leq \frac{1}{\sqrt{A}\cos(\varphi_{GU})},
\end{equation}
and
\begin{equation}\label{543}
\mu(\mathcal{P}_{G P_U(G)^\perp})=\frac{1}{\cos(\varphi_{GU})}
\end{equation}
\begin{proof}
Equation \ttt{543}  follows from \textnormal{(\ref{4.5})}. 
From the definition of $\cos(\varphi_{GU})$ we know that
\begin{equation} \label{11}
\|g\|\cos(\varphi_{GU})\leq \|P_Ug\| \quad \text{for all }g\in G
\end{equation}
Furthermore, from the Cauchy-Schwarz inequality, and \textnormal{(\ref{figs})}, it follows that for 
$\tilde{f}=\mathcal{P}_{G P_U(G)^\perp}f$
\begin{equation}
\langle P_U\tilde{f},\tilde{f}\rangle = \langle P_Uf,\tilde{f}\rangle 
\leq \langle P_Uf,f\rangle^{\frac{1}{2}} \langle P_U\tilde{f},\tilde{f}\rangle^{\frac{1}{2}}.
\end{equation}
This yields
\begin{equation} \label{22}
\|P_U\tilde{f}\|\leq \|P_Uf\|.
\end{equation}
From the definition of a frame sequence,
\begin{equation} \label{fsu}
A\|{u}\|^2\leq \|\U^*u\|^2\leq B\|u\|^2,\quad \text{for }u\in U.
\end{equation}
From \ttt{fsu}, it follows that
\begin{equation} \label{33}
\sqrt{A}~\|P_Uf\|\leq \|\U^*f\|_{l^2}.
\end{equation}
We combine \textnormal{(\ref{11}),(\ref{22})} and \textnormal{(\ref{33})} and obtain
\begin{equation}
\|\tilde{f}\|\cos(\varphi_{GU})\leq \|P_U\tilde{f}\|\leq \|P_Uf\|\leq \frac{1}{\sqrt{A}}~\|\U^*f\|.
\end{equation}
The lower bound of $\textnormal{(\ref{21})}$ follows from 
\begin{equation}
 \|\tilde{f}\|\geq \|P_U\tilde{f}\|\geq \frac{1}{\sqrt{B}}~\|\U^*P_U\tilde{f}\|=\frac{1}{\sqrt{B}}~\|\U^*\tilde{f}~\|,
\end{equation}
where we use \ttt{fsu} for the second inequality.
\end{proof}
\end{thm}

\section{Comparison with generalized sampling}\label{comparison}

We review some important properties of the oblique projection $\PP_{G~
  \mathcal{S}(G)^\perp}$ that was  introduced in
\cite{1,adcock,adhapo12}, and we 
compare it with the oblique projection $\mathcal{P}_{G
  P_U(G)^\perp}$. 

\begin{defi}
Let $cos(\varphi_{GU})>0$. We call the oblique projection $\PP_{G~ \mathcal{S}(G)^\perp}$ generalized sampling.
\end{defi}

We recall that 
calculating the coefficients for the oblique projection $\mathcal{P}_{GP_U(G)^\perp}$ amounts to computing the minimal norm element of
\ttt{leastsquare}. By contrast, for calculating
the coefficients of generalized sampling, we have to calculate the minimal norm element of \ttt{leastsquare1}.
Thus  generalized sampling does not require the additional calculation of $(\U^*\U)^\frac{\dagger}{2}$.

While in general the oblique projections
$\mathcal{P}_{GP_U(G)^\perp}$ and $\PP_{G ~\mathcal{S}(G)^\perp}$ are
rather different, they coincide in several situations. 

The following Lemma can be found in \cite[Lemma 3.7]{adhapo12}
\begin{lem} \label{summ}
Let $G$ and $U$ be finite dimensional subspaces of $H$ with $\dim(G)=dim(U)$.
If $\cos(\varphi_{GU})>0$,
then $G \oplus U^\perp = H$. 
\end{lem}

The following Lemma can be found in \cite[Lemma 4.1]{adhapo12}.
\begin{lem} \label{gencon}
Let $G$ and $U$ be finite dimensional with $\dim(G) = \dim(U)$.
If $\cos(\varphi_{GU})>0$, then 
\begin{equation}\nonumber
\PP_{G ~\mathcal{S}(G)^\perp} = \PP_{G U^\perp},
\end{equation}
i.e. the generalized sampling is exactly the consistent reconstruction.
\end{lem}

Similarly to Lemma \ref{gencon}, we prove the following lemma.
\begin{lem} \label{gencon1}
Let $G$ and $U$ be finite dimensional with $\dim(G) = \dim(U)$.
If $\cos(\varphi_{GU})>0$, then 
\begin{equation}\nonumber
\PP_{G ~P_U(G)^\perp} = \PP_{G U^\perp},
\end{equation}
i.e. the frame independent generalized sampling is exactly the consistent reconstruction.
\begin{proof}
From Lemma \ref{summ}, we infer that $\PP_{G U^\perp}: H\rightarrow G$ is a well defined and bounded mapping.
Clearly $P_U(G)\subset U$. If $P_U g=0$ for some  $g\in G$, then $g=0$, because otherwise
$\cos(\varphi_{GU})=0$. From the injectivity, it follows that $P_U(G)$ is a $n$-dimensional subspace of $U$.
Since the dimension of $U$ is also $n$, we deduce $P_U(G) = U$.
\end{proof}
\end{lem}


\begin{lem}\label{cointight}
Let $\{u_j\}_{j=1,\dots,n}$ and $\{g_k\}_{k=1,\dots,m}$ be finite sequences and $H$ and
let $\cos(\varphi_{GU})>0$. If  $\{u_j\}_{j=1,\dots,m}$ is a tight frame sequence, then 
$\PP_{G ~\mathcal{S}(G)^\perp} = \PP_{G ~P_U(G)^\perp}$.
\begin{proof}
Let $\mathcal{S}$ denote the frame operator of $\{u_j\}_{j=1,\dots,m}$ and $A$ the frame bound. Since $\{u_j\}_{j=1,\dots,m}$ 
is a tight frame sequence for $U$, we have 
\begin{equation}\label{108}
\mathcal{S} u=Au \quad \text{for }u\in U.
\end{equation}
Every element $g\in G$ can be decomposed into $g= g_U+g_{U^\perp}$ with $g_U:= P_U g\in G$ and 
$g_{U^\perp}:= g-P_U g\in U^\perp$. 
From  \eqref{108}  we deduce that $\mathcal{S}g = A g_U$. Consequently,
$A P_U g = \mathcal{S}g$ and $\mathcal{S}(G) = P_U(G)$.
\end{proof}
\end{lem}

Lemma \ref{gencon} in combination with Lemma \ref{gencon1} show that if $\dim(G) = \dim(U)$, then generalized sampling 
and frame independent generalized sampling coincide. Lemma \ref{cointight} shows that they coincide, whenever
$\{u_j\}_{j=1,\dots,m}$ is a tight frame sequence. This is important, because in this case the calculation of 
$(\U^*\U)^\frac{\dagger}{2}$ is not necessary.

In terms of $\cos(\varphi_{GU})$, a bound for the quasi-optimality constant $\mu(\PP_{G ~\mathcal{S}(G)^\perp})$,
is stated in the following 
lemma, see \cite[Corollary 4.3]{adhapo12}.
\begin{lem}\label{abh}
If $\cos(\varphi_{GU})>0$, then 
\begin{equation} \label{mu1}
 1\leq \mu(\PP_{G ~\mathcal{S}(G)^\perp})\leq \frac{\sqrt{B}}{\sqrt{A} \cos(\varphi_{GU})}.
\end{equation}
\end{lem}

In contrast to $\PP_{G~ \mathcal{S}(G)^\perp}$, the nullspace of $\mathcal{P}_{G P_U(G)^\perp}$ does not depend 
on the frame $\{u_j\}_{j\in \mathbb{N}}$, and consequently $\mu(\mathcal{P}_{G P_U(G)^\perp})$ is 
independent of the frame $\{u_j\}_{j\in \mathbb{N}}$.

\begin{figure}[h!]
\centering
\includegraphics[width=0.755\textwidth]{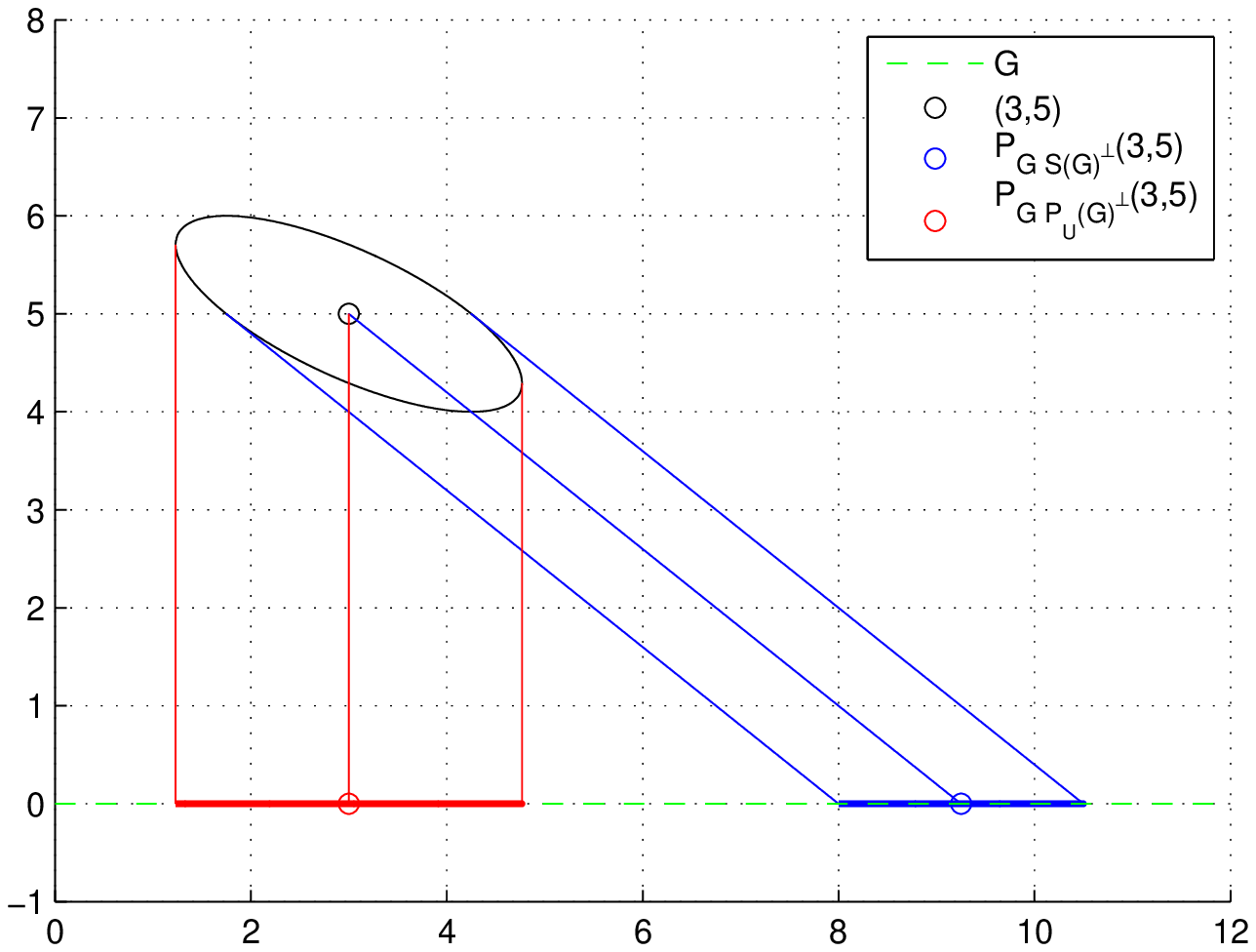}
\caption{$\mathcal{P}_{G P_U(G)^\perp}$ versus $\PP_{G~ \mathcal{S}(G)^\perp}$ for $u_1 = (0,1)$ and
$u_2 = (\frac{4}{5},1)$}
\label{fig_svdcompare}
\end{figure}

The following examples illustrate the difference between $\mathcal{P}_{G P_U(G)^\perp}$ and $\PP_{G~ \mathcal{S}(G)^\perp}$.
Let $H=\mathbb{R}^2$, $u_1 = (0,1)$, $u_2 = (\frac{4}{5},1)$, $g = (1,0)$ and $p = (3,5)$. With that choice,
$G$ is the $x$-axis, $U$ is the whole space $\mathbb{R}^2$ and, consequently, $\mathcal{P}_{G P_U(G)^\perp} = P_G$, the orthogonal 
projection onto $G$. The ellipse in Figure \ref{fig_svdcompare} is the set $E = \{x\in \mathbb{R}^2:\|\U^*(p-x)\|\leq 1\}$. 
We recall, that
\begin{equation}
 \eta(\mathcal{P}_{G P_U(G)}) = \|Q_{|\R(\U^*)}\|,
\end{equation}
where $Q$ is defined by \ttt{Q}.
We observe that 
\begin{align}\label{this}
 \|Q_{|\R(\U^*)}\| &= \underset{\|\U^*f\|=1}{\textnormal{sup}}~\|Q\U^*f\| = 
 \underset{\|\U^*f\|= 1}{\textnormal{sup}}~\|Q\U^*p-Q\U^*(p-f)\|  \\
 &=\underset{\|\U^*(p-x)\|= 1}{\textnormal{sup}}~\|Q\U^*p-Q\U^*(x)\| \\
 &= \underset{\|\U^*(p-x)\|= 1}{\textnormal{sup}}~\|\mathcal{P}_{G P_U(G)^\perp}p-\mathcal{P}_{G P_U(G)^\perp}x\|,
\end{align}
which shows that half of the length of $\mathcal{P}_{G P_U(G)^\perp}(E)$ (red bold segment on the $x$-axis) is $\eta(\mathcal{P}_{G P_U(G)^\perp})$.
Similarly it is shown that half of the length of $\PP_{G~ \mathcal{S}(G)^\perp}(E)$ (blue bold segment on the $x$-axis) is $\eta(\PP_{G~ \mathcal{S}(G)^\perp})$.

 The length of $\mathcal{P}_{G P_U(G)^\perp}(E)$ is greater than the length of 
$\PP_{G ~\mathcal{S}(G)^\perp}(E)$, which shows
that
\begin{equation}
\eta(\mathcal{P}_{G P_U(G)^\perp}) = \|Q_{|\R(\U^*)}\|>\|{Q_1}_{|\R(\U^*)}\| = \eta(\PP_{G ~\mathcal{S}(G)^\perp}). 
\end{equation}
The mapping $\mathcal{P}_{G P_U(G)^\perp}$ is closer to the orthogonal projection $P_U$ than $\PP_{G~ \mathcal{S}(G)^\perp}$ 
(in fact in this example $\mathcal{P}_{G P_U(G)^\perp} = P_U$), which shows that $\mu(\mathcal{P}_{G P_U(G)^\perp})<
\mu(\PP_{G~ \mathcal{S}(G)^\perp})$. 
In this example $G\subset U$, and consequently $\cos(\varphi_{GU}) = 1$. 
\begin{table}[h]
  \centering
  \begin{tabular}{|l|p{1.5cm}|p{1.5cm}|}
    \hline
    \rule[-2mm]{0mm}{6mm}{\em }  & {\em  $\eta$}
            & {\em $\mu$} \\ \hline \hline
    \rule[-1mm]{0mm}{5mm}  $\mathcal{P}_{G P_U(G)^\perp}$   
             & $1.77$ & $1$ \\ \hline
    \rule[-1mm]{0mm}{5mm}  $\PP_{G ~\mathcal{S}(G)^\perp}$   
            &  $1.25$ & 1.6\\ 
    \hline
  \end{tabular}
\caption{The quantities $\eta$ and $\mu$ of $\mathcal{P}_{G P_U(G)^\perp}$ and $\PP_{G~ \mathcal{S}(G)^\perp}$ for $u_1 = (0,1)$ and
$u_2 = (\frac{4}{5},1)$. }
\label{vergleich} 
\end{table}

Since $\cos(\varphi_{GU})=1$, Lemma \ref{abh} implies that the large $\mu(\PP_{G ~\mathcal{S}(G)^\perp})$ is only possible
due to the large value of $\frac{\sqrt{B}}{\sqrt{A}}$. A large value of $\frac{\sqrt{B}}{\sqrt{A}}$ need not cause a big 
$\mu(\PP_{G ~\mathcal{S}(G)^\perp})$. This is illustrated in the following example. Let $H=\mathbb{R}^2$, $u_1 = (1,0)$, $u_2 = (1,\frac{4}{5})$, $g = (1,0)$ and $p = (3,5)$.
\begin{figure}[h!]
\centering
\includegraphics[width=0.755\textwidth]{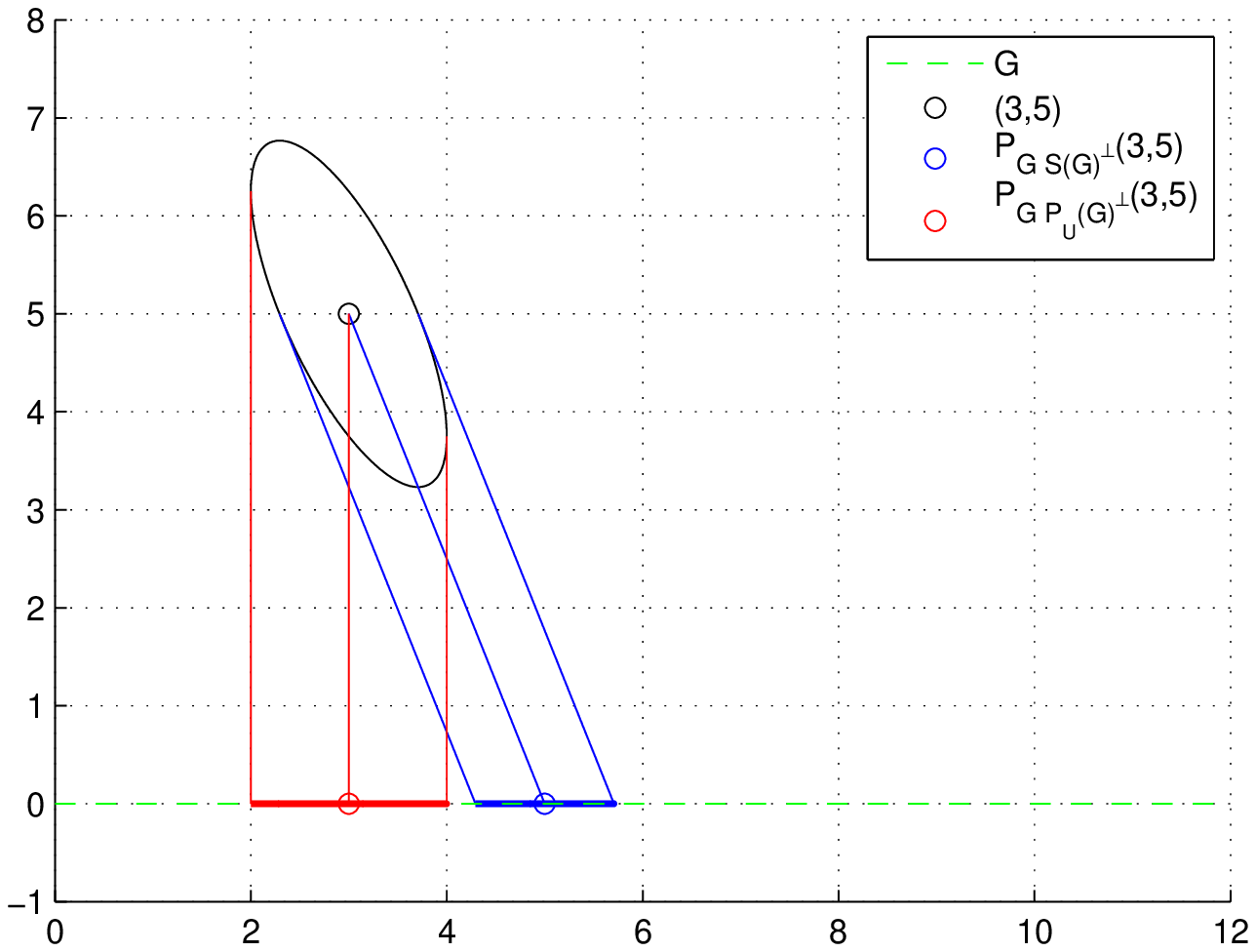}
\caption{$\mathcal{P}_{G P_U(G)^\perp}$ versus $\PP_{G~ \mathcal{S}(G)^\perp}$ for $u_1 = (0,1)$ and
$u_2 = (\frac{4}{5},1)$}
\label{fig_svdcompare1}
\end{figure}
Like in the previous example $\cos(\varphi_{GU}) = 1$. 

\begin{table}[h]
  \centering
  \begin{tabular}{|l|p{1.5cm}|p{1.5cm}|}
    \hline
    \rule[-2mm]{0mm}{6mm}{\em }  & {\em  $\eta$}
            & {\em $\mu$} \\ \hline \hline
    \rule[-1mm]{0mm}{5mm}  $\mathcal{P}_{G P_U(G)^\perp}$   
             & $1$ & $1$ \\ \hline
    \rule[-1mm]{0mm}{5mm}  $\PP_{G ~\mathcal{S}(G)^\perp}$   
            &  $0.71$ & $1.08$\\ 
    \hline
  \end{tabular}
\caption{The quantities $\eta$ and $\mu$ of $\mathcal{P}_{G P_U(G)^\perp}$ and $\PP_{G~ \mathcal{S}(G)^\perp}$ for $u_1 = (0,1)$ and
$u_2 = (\frac{4}{5},1)$.}
\label{vergleich1} 
\end{table}
We observe that in this example the values $\mu(\PP_{G ~\mathcal{S}(G)^\perp})$, 
$\eta(\PP_{G ~\mathcal{S}(G)^\perp})$ and $\eta(\PP_{G P_U(G)^\perp})$ are all smaller than in the previous example. This can
be explained by the fact that $u_1$ and $u_2$ are closer to $G$, and consequently the major axis of the ellipse rotates
into the direction of the $y$ axis.

\end{document}